\documentclass{amsart}

\usepackage{amssymb,latexsym,amsfonts,amsmath}
\usepackage{graphicx,diagrams}
\usepackage{tikz}
\usetikzlibrary{automata}
\usetikzlibrary{shapes}

\newtheorem{theorem}{Theorem}[section]

\newtheorem{problem}[theorem]{Problem}

\newtheorem{proposition}[theorem]{Proposition}
\newtheorem{corollary}[theorem]{Corollary}

\theoremstyle{definition}
\newtheorem{definition}[theorem]{Definition}

\theoremstyle{remark}

\numberwithin{equation}{section}

\newcommand{\R}{{\mathbb{R}}}

\newcommand{\N}{{\mathbb{N}}}

\newcommand{\eg}{{\it e.g. }}

\newcommand{\x}{{\mathbf{x}}}

\newcommand{\Po}{\mathbf{P}}

\begin{document}
\begin{abstract}
The objective of this paper is to solve the controller synthesis problem for bisimulation equivalence in a wide variety of scenarios including discrete-event systems, nonlinear control systems, behavioral systems, hybrid systems and many others. This will be accomplished by showing that the arguments underlying proofs of existence and methods for the construction of controllers are extraneous to the particular class of systems being considered and thus can be presented in greater generality.
\end{abstract}

\title[Controller synthesis for bisimulation equivalence]{Controller synthesis for bisimulation equivalence}
\thanks{This research was partially supported by the National Science Foundation CAREER award 0717188.}

\author[Paulo Tabuada]{Paulo Tabuada}
\address{UCLA Electrical Engineering Department\\
66-147F Engineering IV Building\\
Los Angeles, CA 90095-1594}
\urladdr{http://www.ee.ucla.edu/$\sim$tabuada}
\email{tabuada@ee.ucla.edu}

\maketitle

\section{Introduction}
The notion of bisimulation, introduced by Park~\cite{Park} and Milner~\cite{Milner} in the context of concurrency theory, has been successfully used  as a mechanism to mitigate the complexity of software verification~\cite{ModelChecking}. Recently, the same notion was shown to be relevant for continuous~\cite{BisimTAC,BisimSCL,Gra07}, switched~\cite{SwitchBisim}, hybrid~\cite{BisimTCS} and abstract state systems~\cite{PSDB05}. What makes bisimulation appealing it the possibility of rendering systems of different ``sizes'' equivalent. Here, size needs to be interpreted differently according to the context. When dealing with systems described by finite models, such as discrete-event systems, size means cardinality of the state set. In the case of continuous control systems, size means dimension of the state-space and in the hybrid case size needs to be interpreted as a combination of cardinality and dimension. 

Bisimulation also plays an important role in system synthesis. One can start with a simple model $S$ of a system and try to design a controller $C$ acting on the plant $P$ so that the resulting system $C\parallel P$ is equivalent to $S$. When equivalence is interpreted as isomorphism, the specification $S$ needs to be as complex\footnote{Since most notions of isomorphism are based on an invertible map between the state sets.} as the designed system $C\parallel P$ and this makes this strategy appealing only for small systems. However, when bisimulation is used as equivalence, we can have a specification $S$ being much simpler than the designed system $C\parallel P$. This observation naturally motivates the following controller synthesis problem:

\begin{problem}
\label{CSBE}
Given a plant $P$ and a specification $S$ does there exist a controller $C$ such that the composition $C\parallel P$ is bisimilar to $S$? If so, how do we construct $C$?
\end{problem}


We will solve Problem~\ref{CSBE} in a variety of different contexts thereby recovering known results and proving new ones. The path towards generality followed in this paper is not based on the choice of a model of system that is general enough to contain all the other models as particular cases. Instead, we will work with all the models at the same time. This will be accomplished through the use of elementary ideas from category theory. By proving the results outside any particular class of systems we are able to distill the crucial requirements leading to the existence and the construction of controllers for bisimulation equivalence. The categorical prerequisites are minimal and all the definitions and constructions will be illustrated throughout the paper with transition systems and nonlinear control systems. 

At the technical level we will use the open maps framework of Joyal and coworkers~\cite{OpenMaps} to reason about bisimulation. This framework had already been used in~\cite{OpenTab04} to shown that the controller synthesis problem is solvable in polynomial time for deterministic transition systems and deterministic timed transition systems thus recovering existing results in the computer science literature, see for example~\cite{MadhuTCS} and the references therein. The results of this paper can be seen as a generalization of~\cite{OpenTab04} to a wider class of systems comprising also nonlinear control systems, behavioral systems and hybrid systems.

\section{Notation}
Given a set $S$ we denote by $S^*$ the set of all finite strings obtained by concatenating elements in $S$. An element $s$ of $S^*$ is therefore given by $s=s_1s_2\hdots s_n$ with $s_i\in S\cup\{\epsilon\}$ for $i=1,\hdots,n$ and where $\epsilon$ satisfies $s\epsilon=\epsilon s=s$ for any $s\in S$. 
The length of a string $s\in S^*$ is denoted by $\vert s\vert$. Given a map $f:A\to B$ we shall use the same letter to denote the extension of $f$ to $f:A^*\to B^*$ defined by: 
$$f(s_1s_2\hdots s_n)=f(s_1)f(s_2)\hdots f(s_n).$$

The identity map on a set $A$ will be denoted by $1_A$. 
When $f:M\to N$ is a smooth map between smooth manifolds, $Tf$ will denote the tangent map $Tf:TM\to TN$ taking tangent vectors $X\in T_xM$ at $x\in M$ to tangent vectors $T_xf\cdot X\in T_{f(x)}N$ at $f(x)\in N$. Here $TM=\cup_{x\in M}T_xM$ denotes the tangent bundle of $M$. Map $f$ is said to be a diffeomorphism if there exists a smooth map $g:N\to M$ satisfying $f\circ g=1_N$ and $g\circ f=1_M$.

\section{Systems in categories}
Recall that a category is a collection of objects, that in this paper will model \emph{systems}, and morphisms relating objects. We shall not recall here the precise definitions\footnote{The interested reader is referred to~\cite{MacLane}.} but rather give some simple examples. If one is interested in linear algebra it is natural to take vector spaces as the objects of study and linear maps as morphisms between these objects. If differential geometry is the subject of investigation, objects would be smooth manifolds and smooth maps could be taken as morphisms resulting in the category {\bf Man}. When only the topological structure is of interest, topological spaces would be the objects of study and continuous maps would serve as morphisms. As a final example we mention {\bf Set}, the category having sets as objects and maps between sets as morphisms. To keep the discussion as concrete as possible we will use two examples to illustrate all the definitions and results throughout the paper. The first considers transition systems as a model\footnote{Other models for discrete-event systems are discussed in Section~\ref{SSec:Discrete}.} for discrete-event systems.

\subsection{Transition systems}
A transition system can be seen as a very elementary model of discrete-event systems having while applicability in the formal verification of software~\cite{ModelChecking}.

\begin{definition}
A transition system $T$ is a tuple $T=(Q,\imath,L,\rTo)$ where:
\begin{itemize}
\item $Q$ is a finite set of states;
\item $\imath\in Q$ is the initial state;
\item $L$ is a finite set of labels;
\item $\rTo\subseteq Q\times L\times Q$ is a transition relation.
\end{itemize}
\end{definition}
An element $(p,l,q)\in \rTo$ will be denoted by the more suggestive notation $p\rTo^l q$. When a transition system is used as a model of software, the software execution is described by the notion of run.

\begin{definition}
A run $r$ of a transition system $T=(Q,\imath,L,\rTo)$ is a string $r\in L^*$ for which there exists another string $s\in Q^*$ satisfying:
\begin{enumerate}
\item $s_1=\imath$;
\item $s_i\rTo^{r_i}s_{i+1}$ with $i=1,\hdots,\vert r\vert$.
\end{enumerate}
\end{definition}
A state $q\in Q$ is said to be reachable in $T$ if there exists a run $r$ such that the associated string $s\in Q^*$ satisfies $s_{i+1}=q$.

One possible category for the study of transition systems, denoted by $\mathbf{Tran}$, consists of transition systems as objects and morphisms defined as follows:
\begin{definition}
A morphism $T_1\rTo^f T_2$ from transition system $T_1=(Q_1,\imath_1,L_1,\rTo_1)$ to transition system $T_2=(Q_2,\imath_2,L_2,\rTo_2)$ consists of a pair of maps $f=(f_Q,f_L)$ with $f_Q:Q_1\to Q_2$ and $f_L:L_1\to L_2$ satisfying:
\begin{enumerate}
\item $f_Q(\imath_1)=\imath_2$;
\item $p_1\rTo^{l_1} q_1$ implies $f_Q(p_1)\rTo^{\,\,\, f_L(l_1)\,\,\,} f_Q(q_1)$.
\end{enumerate}
\end{definition}
Other notions of morphism are possible, \eg~\cite{Winskel}, but this one will suffice for our purposes. Note that a morphism from $T_1$ to $T_2$ is guaranteed to take runs of $T_1$ into runs of $T_2$.

\begin{proposition}[Adapted from~\cite{Winskel}]
Let $T_1\rTo^{f} T_2$ be a morphism in {\bf Tran}. Then, for every run $r$ of $T_1$, $f_L(r)$ is a run of $T_2$.
\end{proposition}

\subsection{Control systems}

Nonlinear control systems provide the other example that will be used throughout the paper.

\begin{definition}
A control system $\Sigma$ is a triple $(U,M,F)$ where $U$ is a smooth manifold describing the input space, $M$ is a smooth manifold describing the state space and $F:M\times U\to TM$ is a smooth map describing the system dynamics.
\end{definition}

Trajectories of control systems are defined as usual.

\begin{definition}
A smooth curve $\mathbf{x}:I\to M$ is said to be a trajectory of a control system $\Sigma=(U,M,F)$ if $I\subseteq\R$ is an open interval containing the origin and there exists a smooth curve $\mathbf{u}:I\to U$ satisfying:
$$\frac{d}{dt}\mathbf{x}(t)=F(\mathbf{x}(t),\mathbf{u}(t)),\qquad t\in I$$
\end{definition}

We will say that a control system $\Sigma$ is observable with respect to a smooth map $f:M\times U\to X$ if for any two trajectories $\mathbf{x}$ and $\mathbf{y}$ of $\Sigma$, $\mathbf{x}\ne\mathbf{y}$ implies $f\circ\mathbf{x}\ne f\circ\mathbf{y}$.

The category\footnote{See also~\cite{Elkin,Quotients}.} of control systems, denoted by {\bf Con}, has control systems for objects and morphisms defined as follows.

\begin{definition}
A morphism $\Sigma_1\rTo^{f}\Sigma_2$ from control system $\Sigma_1=(U_1,M_1,F_1)$ to control system $\Sigma_2=(U_2,M_2,F_2)$ consists of a pair of smooth maps $f=(f_M,f_U)$ with $f_M:M_1\to M_2$ and $f_U:M_1\times U_1\to U_2$ satisfying:
\begin{equation}
\label{AffMorphism}
T_xf_M\cdot F_1(x,u)=F_2\big(f_M(x),f_U(x,u)\big)
\end{equation}
\end{definition}

As was the case in {\bf Tran}, morphisms in {\bf Con} transform trajectories into trajectories:

\begin{proposition}[Adapted from~\cite{PLS00}]
Let $\Sigma_1\rTo^{f} \Sigma_2$ be a morphism in {\bf Con}. Then, for every trajectory $\mathbf{x}$ of $\Sigma_1$, $f_M\circ \mathbf{x}$ is a trajectory of $\Sigma_2$.
\end{proposition}

\section{Bisimulation and open maps}
\label{Sec:Bisimulation}
In this section we quickly review the open maps framework introduced by Joyal and co-workers in~\cite{OpenMaps} and apply it to {\bf Tran} and {\bf Con}.

\subsection{General theory}
We consider a category ${\mathbf S}$ of systems with morphisms $X\rTo^{f}Y$ describing how system $Y$ simulates system $X$. In this framework, the notion of bisimulation is introduced by resorting to the notion of path. We thus consider a subcategory ${\mathbf P}$ of ${\mathbf S}$ of path objects whose morphisms describe how paths objects can be extended. Bisimulation is now described through morphisms possessing a special path lifting property:

\begin{definition}
\label{Def:POpen}
A morphism $X\rTo^{f}Y$ is said to be $\Po$-open if given the following commutative diagram:
\begin{diagram}
C & \rTo^{c} & X\\
\dTo^{e} &&\dTo>{f}\\
D & \rTo^{d} & Y
\end{diagram}
where $C$ and $D$ are path objects, there exists a diagonal morphism $D\rTo^{r} X$ making the following diagram commutative:
\begin{diagram}
C & \rTo^{c} & X\\
\dTo^{e} &\ruTo^{r}&\dTo>{f}\\
D & \rTo^{d} & Y
\end{diagram}
that is, $c=r\circ e$ and $d=f\circ r$.
\end{definition}

\subsection{Examples}
\subsubsection{Transition systems}
\label{SSec:TransitionBisim}
The notion of bisimulation was introduced by Park~\cite{Park} and Milner~\cite{Milner} in the context of transition systems as follows:

\begin{definition}
Let $T_1$ and $T_2$ be transition systems with the same label set $L$. A relation $R\subseteq Q_1\times Q_2$ with $(\imath_1,\imath_2)\in R$ is said to be a simulation relation from $T_1$ to $T_2$ if  $(p_1,p_2)\in R$ implies:
\begin{enumerate}
\item $p_1\rTo^l q_1$ in $T_1$ implies existence of $p_2\rTo^l q_2$ in $T_2$ with $(q_1,q_2)\in R$. 
\label{RI}
\end{enumerate}
A relation $R\subseteq Q_1\times Q_2$ is said to be a bisimulation relation between $T_1$ and $T_2$ if $(p_1,p_2)\in R$ implies in addition to~(\ref{RI}):
\begin{enumerate}
\setcounter{enumi}{1}
\item $p_2\rTo^l q_2$ in $T_2$ implies existence of $p_1\rTo^l q_1$ in $T_1$ with $(q_1,q_2)\in R$.
\end{enumerate}

Transition systems $T_1$ and $T_2$ are said to be bisimilar if there exists a bisimulation relation between them.
\end{definition}

According to this definition, transitions in $T_1$ must be matched by transitions in $T_2$ with the same label and, conversely, transitions in $T_2$ must be matched by transitions  in $T_1$ also with the same label. To capture this requirement on the labels, using the open maps framework, we fix a set of labels $L$ and let {\bf S} be the subcategory {\bf Tran}$_L$ of {\bf Tran} consisting of transition systems with label set $L$ and morphisms $f:T_1\to T_2$ satisfying $f_L=1_L$. For the path subcategory {\bf P} we take the full\footnote{A category {\bf D} is a full subcategory of a category {\bf C} when any object of {\bf D} is also an object of {\bf C} and for any two objects $X$ and $Y$ in {\bf D}, if $X\rTo^f Y$ is a morphism in {\bf C} then it is also a morphism in {\bf D}.} subcategory of {\bf S} defined by objects of the form:
\begin{equation}
\label{Path}
q_1\rTo{l_1}q_2\rTo^{l_2}q_3\rTo^{l_3}\hdots\rTo^{l_{n-1}}q_n
\end{equation}
with $q_1=\imath$ and $q_i\ne q_j$ for $i\ne j$. Note that any morphism $T\rTo^f T_1$ from a path object $T$ describes a run $l_1l_2\hdots l_{n-1}$ of $T_1$ through the sequence of transitions $f_Q(q_i)\rTo^{f_L(l_i)=l_i}f_Q(q_{i+1})$ in $T_1$. Conversely, every run of $T_1$ can be described by a morphism from a path object into $T_1$.

With this choice for {\bf S} and {\bf P} we recover Park~\cite{Park} and Milner's~\cite{Milner} notion of bisimulation through a diagram of $\Po$-open maps.

\begin{theorem}[\cite{OpenMaps}]
Let $T_1$ and $T_2$ be objects in {\bf S}. $T_1$ is bisimilar to $T_2$ iff there exists a diagram:
\begin{equation}
\label{Span}
T_1\lTo^{\alpha}T\rTo^{\beta}T_2
\end{equation}
where $\alpha$ and $\beta$ are $\Po$-open morphisms.
\end{theorem}

The intuition behind the diagram~(\ref{Span}) can be understood by noting that a diagram $C\lTo^{f}B\rTo^{g}D$ in {\bf Set} defines a relation $R\subseteq C\times D$ by $(c,d)\in R$ if there is a $b\in B$ such that $f(b)=c$ and $g(b)=d$. Conversely, given a relation $R\subseteq C\times D$ we can always construct a diagram $C\lTo^f R\rTo^g D$ where $f=\pi_C\circ i$ and $g=\pi_D\circ i$ with $i:R\to C\times D$ being the natural inclusion of $R$ in $C\times D$, and $\pi_C:C\times D\to C$ and $\pi_D:C\times D\to D$ the natural projections. The diagram~(\ref{Span}) is then simply defining the relation $R\subseteq Q_1\times Q_2$ with $(q_1,q_2)\in Q_1\times Q_2$ if there exists a $q\in Q$ such that $\alpha_Q(q)=q_1$ and $\beta_Q(q)=q_2$. Since $\alpha$ is $\Po$-open, transitions in $T_1$ can be lifted, as described in Definition~\ref{Def:POpen}, to $T$ and then mapped to $T_2$ through the morphism $\beta$. We thus see that $\Po$-openness of $\alpha$ ensures that $R$ is a simulation relation from $T_1$ to $T_2$. Moreover, as $\beta$ is also $\Po$-open, transitions in $T_2$ can also be matched by transitions in $T_1$ thus making $R$ a bisimulation.

\subsubsection{Control systems}
The notion of bisimulation was recently studied in the context of nonlinear control systems~\cite{BisimTAC,BisimSCL,BisimTCS}. In this paper we formalize bisimulation for control systems as follows:

\begin{definition}[Adapted from~\cite{BisimSCL,BisimTCS}]
\label{BisimControl}
Let $\Sigma_1=(U_1,M_1,F_1)$ and $\Sigma_2=(U_2,M_2,F_2)$ be control 
systems and let $R\subseteq M_1\times M_2$ be a submanifold of $M_1\times M_2$ for which the natural projection maps $\pi_1:R\to M_1$ and $\pi_2:R\to M_2$ are surjective submersions. Relation $R\subseteq M_1\times M_2$ is said to be a simulation relation 
from $M_1$ to $M_2$ if $(x_1,x_2)\in R$ implies:
\begin{enumerate}
\item for any trajectory $\x_1:I\to M_1$ of $\Sigma_1$ with $\x_1(0)=x_1$ there exists a  trajectory $\x_2: I\to M_2$ of $\Sigma_{2}$ with $\x_2(0) = x_2$ such that $(\x_1(t), \x_2(t))\in R$ for every $t\in I\cap\R_0^+$.
\label{RelI}
\end{enumerate}

A relation $R\subseteq M_1\times M_2$ is said to be a bisimulation relation between $\Sigma_1$ and $\Sigma_2$ if $(x_1,x_2)\in R$ implies in addition to~(\ref{RelI}):

\begin{enumerate}
\setcounter{enumi}{1}
\item for any trajectory $\x_2:I\to M_2$ of $\Sigma_2$ with $\x_2(0)=x_2$ there exists a trajectory $\x_1: I\to M_1$ of $\Sigma_{1}$ with $\x_1(0) = x_1$ such that $(\x_1(t), \x_2(t))\in R$ for every $t\in I\cap\R_0^+$.
\end{enumerate}
\end{definition}

When control systems $\Sigma_1$ and $\Sigma_2$ are equipped with observation maps $h_1:M_1\to O$ and $h_2:M_2\to O$, respectively, the above notion can be strengthened by requiring that states $(x_1,x_2)\in R$ also satisfy $h_1(x_1)=h_2(x_2)$. This is the approach taken in~\cite{BisimTAC} which can also be captured in the proposed framework by defining a category of control systems equipped with observation maps.

Definition~\ref{BisimControl} requires $R$ to be a manifold and the projection maps $\pi_i:R\to M_i$ to be surjective submersions. Although the notion of bisimulation still makes sense without these technical requirements, they are used to guarantee\footnote{The open maps approach requires a category with finite pullbacks. {\bf Con} is based on {\bf Man} since the state and input spaces are manifolds and in {\bf Man} pullbacks do not always exist. This can be remedied by using surjective submersions for which pullbacks are guaranteed to exist.} that bisimulation is a notion of equivalence in~{\bf Con} as discussed in~\cite{BisimTCS}.
 
In order to describe bisimulations in {\bf Con} through open maps we take {\bf S}$=${\bf Con} and consider the full subcategory of {\bf Con} defined by objects of the form $\Sigma=(\{*\},I,F)$ where $\{*\}$ is a set with a single element $*$, $I\subseteq\R$ is an open interval containing the origin and $F$ is defined by $F(t,*)=F(t)=1$ for any $t\in I$. Intuitively, $\Sigma$ describes time modeled as a control system. A morphism $\Sigma\rTo^f \Sigma_1$ from a path object $\Sigma$ is described by a pair of smooth maps $f_M:I\to M$ and $f_U:I\to U$ satisfying~(\ref{AffMorphism}):
$$\frac{d}{dt} f_M(t)=T_tf_M\cdot 1=T_tf_M\cdot F(t)=F_1(f_M(t),f_U(t))$$
We thus see that a morphism $\Sigma\rTo^f \Sigma_1$ from a path object describes a trajectory $f_M:I\to M$ of $\Sigma_1$ induced by the input curve $f_U:I\to U$. Conversely, every trajectory of $\Sigma_1$ can be seen as a morphism from a path object into $\Sigma_1$.

With this choice for {\bf S} and {\bf P} we have the following result:

\begin{theorem}[\cite{BisimTCS}]
Let $\Sigma_1$ and $\Sigma_2$ be objects in {\bf Con}. $\Sigma_1$ is bisimilar to $\Sigma_2$ iff there exists a diagram:
$$\Sigma_1\lTo^{\alpha}\Sigma\rTo^{\beta}\Sigma_2$$
where $\alpha$ and $\beta$ are $\Po$-open morphisms with $\alpha_M$ and $\beta_M$ are surjective submersions.
\end{theorem}

\section{Composition as a pullback}
Before addressing problems of control we need one last ingredient: composition of systems. Although composition assumes very different forms for different classes of systems we can obtain a unified description by resorting to the notion of pullback. The use of pullbacks to describe system composition has been used several times before, \eg~\cite{BBCPP03,ASS06}.

\subsection{General theory}

\begin{definition}
\label{Def:Pullback}
The pullback of two morphisms $X\rTo^{x_a} A$ and $Y\rTo^{y_a} A$ in a category is a pair of morphisms $Z\rTo^\alpha X$ and $Z\rTo^\beta Y$ satisfying $x_a\circ \alpha=y_a\circ \beta$ and such that for any other pair of morphisms $Z'\rTo^{\alpha'} X$ and $Z'\rTo^{\beta'} Y$ satisfying $x_a\circ \alpha'=y_a\circ \beta'$ there exists a unique morphism $Z'\rTo^\gamma Z$ making the following diagram commutative:
\begin{diagram}[small]
&&Z'\\
&\ldTo(2,4)^{\alpha'}&\dTo^{\gamma}&\rdTo(2,4)^{\beta'}\\
&&Z\\
&\ldTo>{\alpha}&&\rdTo<{\beta}\\
X&&&&Y\\
&\rdTo^{x_a}&&\ldTo^{y_a}\\
&&A
\end{diagram}
\end{definition}

The pullback of $X\rTo^{x_a} A$ and $Y\rTo^{y_a} A$ is denoted by $X\times_A Y$. When the pullback of any two morphisms $X\rTo^{x_a} A$ and $Y\rTo^{y_a} A$ in a category {\bf S} exists we say that {\bf S} has binary pullbacks. As with many other definitions in category theory, pullbacks are uniquely defined up to isomorphism. This means that any two objects $Z_1$ and $Z_2$ satisfying the above definition are necessarily isomorphic in the sense that there exist morphisms $f:Z_1\to Z_2$ and $g:Z_2\to Z_1$ satisfying $f\circ g=1_{Z_2}$ and $g\circ f=1_{Z_1}$.

Pullbacks $X\times_A Y$ in {\bf Set} can be constructed by first computing the Cartesian product $X\times Y$ and then selecting the elements of $(x,y)\in X\times Y$ satisfying the equality $x_a(x)=y_a(y)$. $X\times_M Y$ is then given by the set $\{(x,y)\in X\times Y\,\,\vert\,\, x_a(x)=y_a(y)\}$ equipped with the maps $\alpha=\pi_X\circ i$ and $\beta=\pi_Y\circ i$ where $i:X\times_M Y\to X\times Y$ is the natural inclusion of $X\times_A Y$ into $X\times Y$, and $\pi_X:X\times Y\to X$ and $\pi_Y:X\times Y\to Y$ are the natural projections. We leave to the reader to verify that $X\times_A Y$ constructed as described above does satisfy Definition~\ref{Def:Pullback}. The same idea underlies the construction of pullbacks in {\bf Tran} and {\bf Con} as described later in this section. The object $A$ serves as a mediator or interface between the objects $X$ and $Y$. By changing $A$, $x_a$ and $y_a$ we can model a wide variety of interconnections between systems. When $A$ is seen as an interface we can regard the morphisms $x_a$ and $y_a$ the description of how the internal state is exposed through the interface. In this way, the pullback $X\times_A Y$ describes the result of interconnecting $X$ to $Y$ through the interface $A$. However, more interesting types of interconnection, such as feedback, can still be modeled by pullbacks as we next describe.

\subsection{Examples}
\subsubsection{Transition systems}
The most frequently used composition of transition systems requires synchronization on common labels or events. 

\begin{definition}
Let $T_1=(Q_1,\imath_1,L,\rTo_1)$ and $T_2=(Q_2,\imath_2,L,\rTo_2)$ be transition systems. The parallel composition of $T_1$ and $T_2$, denoted by $T_1\parallel T_2$, is the transition system $T_1\parallel T_2=(Q_{12},\imath_{12},L_{12},\rTo_{12})$ defined by:
\begin{itemize}
\item $Q_{12}=Q_1\times Q_2$;
\item $\imath_{12}=(\imath_1,\imath_2)$;
\item $L_{12}=L$;
\item $(p_1,p_2)\rTo_{12}^l (q_1,q_2)$ in $T_1\parallel T_2$ if $p_1\rTo_{1}^l q_1$ in $T_1$ and $p_2\rTo_{2}^l q_2$ in $T_2$.
\end{itemize}
\end{definition}

In order to model $T_1\parallel T_2$ as a pullback in {\bf Tran}$_L$ we first note that given  $T_1\rTo^{t_{1a}}T_A$ and $T_2\rTo^{t_{2a}}T_A$ we can construct $T_1\times_{T_A}T_2$ by first constructing the state set as: 
$$\big\{(q_1,q_2)\in Q_1\times Q_2\,\,\vert\,\, t_{1aQ}(q_1)=t_{2aQ}(q_2)\big\}$$ 
and then constructing the transition relation as:
$$\big\{((p_1,p_2),l,(q_1,q_2))\in Q\times L\times Q\,\,\vert\,\, (t_{1aQ}(p_1),l,t_{1aQ}(q_1))=(t_{2aQ}(p_2),l,t_{2aQ}(q_2)) \big\}$$

Using this insight we define the transition system $T_A=(Q_A,\imath_A,L_A,\rTo_A)$: 
$$Q_A=\{*\},\qquad \imath_A=*,\qquad L_A=L,\qquad \rTo_A=\bigcup_{l\in L}\big\{(*,l,*)\big\}$$
and note that for any transition system $T=(Q,\imath,L,\rTo)$ there exists a morphism $T\rTo^{t_a} T_A$ defined by $t_{aQ}(q)=*$ for every $q\in Q$ and $t_{aL}=1_L$. We now have the following description of $T_1\parallel T_2$:

\begin{proposition}
Let $T_1$ and $T_2$ be transition systems with label set $L$. Then, the parallel composition $T_1\parallel T_2$ is the pullback of $T_1\rTo^{t_{1a}}T_A$ and $T_2\rTo^{t_{2a}}T_A$ in {\bf Tran}$_L$.
\end{proposition}

\begin{proof}[Proof sketch]
$T_1\times_{T_A} T_2$ is equipped with morphisms $T_1\times_{T_A} T_2\rTo^\alpha T_1$ and $T_1\times_{T_A} T_2\rTo^\beta T_2$ defined by $\alpha_Q=\pi_{Q_1}$, $\alpha_L=1_L$, $\beta_Q=\pi_{Q_2}$ and $\beta_L=1_L$ where $\pi_{Q_i}:Q_1\times Q_2\to Q_i$ are the natural projections. It is not difficult to verify that $t_{1a}\circ\alpha=t_{1b}\circ\beta$. Let now $T$ be a transition system equipped with morphisms $T\rTo^{\alpha'} T_1$ and $T\rTo^{\beta'} T_2$ satisfying $t_{1a}\circ\alpha'=t_{2a}\circ\beta'$. We now show existence of a unique morphism $T\rTo^\gamma T_1\times_{T_A}T_2$ satisfying $\alpha\circ \gamma=\alpha'$ and $\beta\circ\gamma=\beta'$. Since in {\bf Tran}$_L$ every morphism $f$ has $f_L=1_L$ we conclude that $\gamma_L=1_L$. Moreover, we define $\gamma_Q$ by $\gamma_Q(q)=\big(\alpha'_Q(q),\beta'_Q(q)\big)$. It then follows that $\alpha_Q\circ \gamma_Q(q)=\alpha_Q\big(\alpha'_Q(q),\beta'_Q(q)\big)=\pi_{Q_1}\big(\alpha'_Q(q),\beta'_Q(q)\big)=\alpha'_Q(q)$. Similarly, $\beta_Q\circ \gamma_Q(q)=\beta_Q\big(\alpha'_Q(q),\beta'_Q(q)\big)=\pi_{Q_2}\big(\alpha'_Q(q),\beta'_Q(q)\big)=\beta'_Q(q)$. Assume now that $\gamma$ is not unique and let $\gamma'$ be another morphism from $T$ to $T_1\times_{T_A}T_2$ satisfying $\alpha\circ \gamma'=\alpha'$ and $\beta\circ\gamma'=\beta'$. Then, $\pi_{Q_1}\circ \gamma'_Q(q)=\alpha_Q\circ\gamma'_Q(q)=\alpha'_Q(q)=\alpha_Q\circ\gamma_Q(q)=\pi_{Q_1}\circ \gamma_Q(q)$ and $\pi_{Q_2}\circ \gamma'_Q(q)=\beta_Q\circ\gamma'_Q(q)=\beta'_Q(q)=\beta_Q\circ\gamma_Q(q)=\pi_{Q_2}\circ \gamma_Q(q)$. Since $\pi_{Q_1}\circ \gamma_Q'=\gamma_Q$ and $\pi_{Q_2}\circ \gamma_Q'=\gamma_Q$ we conclude that $\gamma'_Q=\gamma_Q$. The equality $\gamma'=\gamma$ now follows from $\gamma'_L=1_L=\gamma_L$.
\end{proof}

The mediator $T_A$ is rather special in that any string $r\in L^*$ is a run of $T_A$. This choice for $T_A$ was designed to guarantee that runs of $T_1\times_{T_A} T_2$ are the intersection of the runs of $T_1$ and $T_2$. Note that a run of $T_1\times_{T_A} T_2$ should be a pair $(r,s)$ where $r$ is a run of $T_1$ and $s$ is a run of $T_2$ that satisfy $r=t_{1aL}(r)=t_{1aL}(s)=s$. We can, therefore, identify these pairs with the runs $r=s\in L^*$ of $T_1$ and $T_2$. The next section will use very different choices for the mediating object since we are no longer interested in the intersection of behaviors but rather on feedback.

\subsubsection{Control systems}
Control systems can be composed in many different ways. In this section we focus our attention on feedback interconnections. The first kind of interconnection describes the effect of applying a feedback control law $u=k(x,v)$ to a control system $F(x,u)$ resulting in the closed loop system described by $F(x,k(x,v))$. Note that as a special case we have control laws $u=k(x)$ resulting in closed loop systems $F(x,k(x))$ which are no longer affected by the input. 

\begin{definition}
\label{Def:ControlInterconnectI}
Let $\Sigma_1=(U_1,M_1,F_1)$ be a control system and let $k_2:M_1\times U_2\to U_1$ be a smooth feedback law. The feedback interconnection between $\Sigma_1$ and $k_2$ is the control system $\Sigma=(U,M,F)$ with $U=U_2$, $M=M_1$ and $F(x_1,u_2)=F_1(x_1,k(x_1,u_2))$ for every $x_1\in M_1$ and $u_2\in U_2$.
\end{definition}

The second kind of interconnection models the effect of dynamic feedback.

\begin{definition}
\label{Def:ControlInterconnectII}
Let $\Sigma_1=(U_1\times V_1,M_1,F_1)$ and $\Sigma_2=(U_2\times V_2,M_2,F_2)$ be control systems. The feedback interconnection between $\Sigma_1$ and $\Sigma_2$, with interconnection maps $\phi_1:M_1\to U_2$ and $\phi_2:M_2\to U_1$, is the control system $\Sigma=(U,M,F)$ with $U=V_1\times V_2$, $M=M_1\times M_2$ and $F(x,u)=\big(F_1(x_1,(\phi_2(x_2),v_1)),F_2(x_2,(\phi_1(x_1),v_2))\big)$ for every $x_1\in M_1$, $x_2\in M_2$, $v_1\in V_1$ and $v_2\in V_2$.
\end{definition}

Feedback interconnections can be seen as pullbacks by properly defining the mediating object $\Sigma_A$ and the morphisms $\Sigma_1\rTo^{\sigma_{1a}}\Sigma_A$ and $\Sigma_2\rTo^{\sigma_{2a}}\Sigma_A$ as shown in the next propositions.

\begin{proposition}
\label{PropCompose}
Let $\Sigma_1=(U_1,M_1,F_1)$ and $\Sigma_2=(U_2,M_2,F_2)$ be two objects in {\bf Con} where $M_2=M_1$ and $F_2(x_2,u_2)=F_1(x_2,k(x_2,u_2))$ for a smooth feedback law $k:M_2\times U_2\to U_1$. The feedback interconnection of $\Sigma_1$ with $k$ is the pullback of $\Sigma_1\rTo^{\sigma_{1a}} \Sigma_A$ and $\Sigma_2\rTo^{\sigma_{2a}} \Sigma_A$ where $\Sigma_A=\Sigma_1$, $\sigma_{1aM}(x_1)=x_1$, $\sigma_{2aM}(x_2)=x_2$, $\sigma_{1aU}(x_1,u_1)=u_1$ and $\sigma_{2aU}(x_2,u_2)=k(x_2,u_2)$ for every $x_1\in M_1$, $x_2\in M_2$, $u_1\in U_1$ and $u_2\in U_2$.
\end{proposition}

\begin{proof}[Proof sketch]
The result follows by noting that the state space of $\Sigma_1\times_{\Sigma_A}\Sigma_2$ is the set of pairs $(x_1,x_2)\in M_1\times M_2$ satisfying $\sigma_{1aM}(x_1)=\sigma_{2aM}(x_2)$. Since $\sigma_{1aM}=1_{M_1}=1_{M_2}=\sigma_{2aM}$ we can identify $M$ with $M_1=M_2$ through $(x,x)\leftrightarrow x$. The input space is the set of pairs $(u_1,u_2)\in U_1\times U_2$ satisfying $\sigma_{1aU}(x_1,u_1)=\sigma_{2aU}(x_2,u_2)$, or equivalently, $u_1=k(x_2,u_2)$. We can thus identify the set of inputs with $U_2$ through $(k(x_2,u_2),u_2)\leftrightarrow u_2$. Finally, $F$ will be given by the restriction of $(F_1(x_1,u_1),F_2(x_2,u_2))$ to $M\times U$ which can be identified with the points $((x,x),((k(x,u),u))\in (M_1\times M_2)\times(U_1\times U_2)$ thus leading to $F(x,u)=F_1(x,k(x,u))$.
\end{proof}

\begin{proposition}
Let $\Sigma_1=(U_1\times V_1,M_1,F_1)$ and $\Sigma_2=(U_2\times V_2,M_2,F_2)$ be two objects in {\bf Con} and consider the object $\Sigma_A=(U_A,M_A,F_A)$ with $U_A=U_1\times U_2$, $M_A=\{*\}$ and $F_A(*,u_a)=0$ for every $u_a\in U_A$. The feedback interconnection of $\Sigma_1$ and $\Sigma_2$, with interconnection maps $\phi_1:M_1\to U_2$ and $\phi_2:M_2\to U_1$, is the pullback of $\Sigma_1\rTo^{\sigma_{1a}} \Sigma_A$ and $\Sigma_2\rTo^{\sigma_{2a}} \Sigma_A$ where $\sigma_{1aM}(x_1)=*$, $\sigma_{2aM}(x_2)=*$, $\sigma_{1aU}\big(x_1,(u_1,v_1)\big)=(u_1,\phi_1(x_1))$ and $\sigma_{2aU}\big(x_2,(u_2,v_2)\big)=(\phi_2(x_2),u_2)$ for every $x_1\in M_1$, $x_2\in M_2$, $u_1\in U_1$, $u_2\in U_2$, $v_1\in V_1$ and $v_2\in V_2$.
\end{proposition}

\begin{proof}[Proof sketch]
Similar to the proof sketch of Proposition~\ref{PropCompose}.
\end{proof}

Note that arbitrary pullbacks do not exist in {\bf Con} since they do not exist in {\bf Man}. However, the above defined pullbacks are guaranteed to exist.

\section{Existence and synthesis of controllers}
\subsection{General theory}
We now consider the control synthesis problem for bisimulation equivalence. We assume that we are given:
\begin{enumerate}
\item a morphism $P\rTo^{p_a} A$ describing the plant $P$ and the mediator $A$ to be used for control;
\item a morphism $S\rTo^{s_a} A$ describing the specification $S$ and how it relates to the mediator $A$.
\end{enumerate}
Based on this data we formulate the notion of controller as follows:

\begin{definition}
\label{BisimCont}
Let $P\rTo^{p_a}A$, $S\rTo^{s_a}A$ and $C\rTo^{c_a}A$ be morphisms in a category {\bf S}. The morphism $C\rTo^{c_a}A$ is a bisimulation controller for plant $P\rTo^{p_a}A$, enforcing specification $S\rTo^{s_a}A$, if there exists a commutative diagram:
\begin{equation}
\label{SqI}
\begin{diagram}[small]
&&Z&&\\
&\ldTo^s&&\rdTo^{cp}\\
S&&&&C\times_A P\\
&\rdTo^{s_a}&&\ldTo^{cp_a}\\
&&A
\end{diagram}
\end{equation}
 in which $s$ and $cp$ are $\Po$-open morphisms.
\end{definition}

The diagram $S\lTo^s Z\rTo^{cp} C\times_A P$ of $\Po$-open morphisms in diagram~(\ref{SqI}) requires the closed loop system $C\times_A P$ to be bisimilar to $S$. Moreover, commutativity of~(\ref{SqI}) imposes the additional requirement that any two states related through the relation defined by the diagram $S\lTo^s Z\rTo^{cp} C\times_A P$ are indistinguishable by the mediator. This is a natural requirement since both the specification $S$ and the controlled system $C\times_A P$ should behave in the same way when composed with other systems through the mediator $A$. 

We now introduce what can be seen as an observability property.

\begin{definition}
Let $X\rTo^f Y$ be a morphism in {\bf S}. We say that $f$ is a $\Po$-faithfull morphism if given the following commutative diagram:
\begin{diagram}
C & \rTo^{c} & X\\
\dTo^{e} &&\dTo>{f}\\
D & \rTo^{d} & Y
\end{diagram}
where $C$ and $D$ are objects in {\bf P}, existence of diagonal morphisms $D\rTo^{r} X$ and $D\rTo^{s} X$ making the following two diagrams commutative:
\begin{diagram}
C & \rTo^{c} & X && C & \rTo^{c} & X\\
\dTo^{e} &\ruTo^{r}&\dTo>{f} && \dTo^{e} &\ruTo^{s}&\dTo>{f}\\
D & \rTo^{d} & Y && D & \rTo^{d} & Y
\end{diagram}
implies $r=s$.
\end{definition}

We postpone until Section~\ref{SSec:Existence-Examples} a discussion of $\Po$-faithfulness in the concrete context of transition systems and control systems. The main contribution of this paper can now be stated as follows.

\begin{theorem}
\label{ExistBisimContI}
Let $P\rTo^{p_a}A$ and $S\rTo^{s_a}A$ be morphisms in a category {\bf S} with binary pullbacks and assume that $P\rTo^{p_a}A$ is $\Po$-faithfull. There exists a bisimulation controller $C\rTo^{c_a}A$ for plant $P\rTo^{p_a}A$ enforcing specification $S\rTo^{s_a}A$ iff there is a commuting diagram
\begin{equation}
\label{SpanI} 
\begin{diagram}[small]
&&Z&&\\
&\ldTo^\gamma &&\rdTo^\delta\\
S&&&&P\\
&\rdTo^{s_a} &&\ldTo^{p_a}\\
&&A&&
\end{diagram}
\end{equation}
with $\gamma$ a $\Po$-open morphism. Furthermore, when a bisimulation controller $C\rTo^{c_a}A$ exists, we can take $C=S$ and $c_a=s_a$.
\end{theorem}

\begin{proof}
Assume that a bisimulation controller $C$ exists. Then, we have a commuting diagram:
\begin{equation}
\begin{diagram}[small]
&&X&&\\
&\ldTo^s&&\rdTo^{cp}\\
S&&&&C\times_A P\\
&\rdTo^{s_a}&&\ldTo^{cp_a}\\
&&A
\end{diagram}
\end{equation}
where $s$ and $cp$ are $\Po$-open. Taking $Z=X$, $\gamma=s$ and $\delta=p\circ cp$, where $p$ is the morphism $C\times_A P\rTo^p P$, we have a commuting diagram as in~(\ref{SpanI}). Clearly, $\gamma$ is $\Po$-open.

Assume now that diagram~(\ref{SpanI}) exists and let us prove that $C=S$ and $c_a=s_a$ is the desired controller. It follows from the definition of $S\times_A P$ the existence of a unique morphism $Z\rTo^\mu S\times_A P$ satisfying $s\circ \mu=\gamma$ and $p\circ \mu=\delta$. The remaining proof consists in showing that $\mu$ is $\Po$-open since in this case the result follows from the commuting diagram:
\begin{equation}
\begin{diagram}[small]
&&Z&&\\
&\ldTo^\gamma&&\rdTo^{\mu}\\
S&&&&S\times_A P\\
&\rdTo^{s_a}&&\ldTo^{sp_a}\\
&&A
\end{diagram}
\end{equation}
where $\gamma$ is $\Po$-open by assumption. Consider the following commutative diagrams:
\begin{equation}
\label{TwoSq}
\begin{diagram}
C & \rTo^{c} & Z && C & \rTo^c & Z\\
\dTo^{e} &&\dTo>{\mu} && \dTo^e && \dTo>{s\circ \mu}\\
D & \rTo^{d} & S\times_A P && D & \rTo^{s\circ d} & S
\end{diagram}
\end{equation}
where $s$ is the morphism $S\times_A P\rTo^s S$. Since $s\circ \mu=\gamma$ and $\gamma$ is $\Po$-open, there exists a diagonal morphism $D\rTo^r Z$ for the right diagram in~(\ref{TwoSq}). We now show that $D\rTo^r Z$ is also the desired diagonal morphism for the left diagram in~(\ref{TwoSq}). We first note that equality $c=r\circ e$ is inherited from the right diagram so that we only need to show that $\mu\circ r=d$. The equality will be proved by noting that it follows from the fact that $S\times_A P$ is a pull-back that any two morphisms $D\rTo^d S\times_A P$ and $D\rTo^{\mu\circ r} S\times_A P$ are necessarily the same when the following two conditions hold:
\begin{eqnarray}
\label{CondI}
s\circ \mu\circ r & = & s\circ d\\
\label{CondII}
p\circ \mu\circ r & = & p\circ d
\end{eqnarray}
Equality~(\ref{CondI}) follows from the right diagram in~(\ref{TwoSq}) and equality~(\ref{CondII}) follows from $p_a$ being uniquely $\Po$-open and the equality $p_a\circ p\circ \mu\circ r=p_a \circ p\circ d$.

\end{proof}

\subsection{Examples}
\label{SSec:Existence-Examples}
\subsubsection{Transition systems}
For the choice of {\bf S} and {\bf P} described in Section~\ref{SSec:TransitionBisim}, $\Po$-faithfulness of a morphism $T\rTo^{t_a} T_A$ is implied by determinism of $T$. Recall that a transition system is deterministic when $p\rTo^{l} q_1$ and $p\rTo^{l} q_2$ imply $q_1=q_2$. Determinism of $T$ guarantees that a run $r\in L^*$ uniquely determines the string $s\in Q^*$ satisfying $s_i\rTo^{l_i} s_{i+1}$ and thus implies $\Po$-faithfulness of  $T\rTo^{t_a} T_A$. Recalling that a diagram $X\lTo^\alpha Z\rTo^{\beta}Y$ with $\alpha$ a $\Po$-open morphism can be seen as a simulation relation from $X$ to $Y$ we have the following corollary of Theorem~\ref{ExistBisimContI}.

\begin{corollary}
\label{Corol}
Let $T_P$ and $T_S$ be transition systems and assume that $T_P$ is deterministic. There exists a transition system $T_C$ making $T_C\parallel T_P$ bisimilar to $T_S$ iff there exists a simulation relation from $T_S$ to $T_P$.
\end{corollary}
Combining this corollary with existing results on the existence and computation of simulation relations~\cite{SimulEff,SimulHH} we immediately conlude that the controller synthesis problem for  deterministic plants can be solved in polynomial time thus recovering the results in~\cite{MadhuTCS}. In section~\ref{SSec:Discrete} we compare this result with existing results for other models of discrete-event systems.

\begin{figure}[h]
\centering
\begin{tikzpicture}[shorten >=1pt,node distance=2cm,auto]
\tikzstyle{every state}=[draw=blue!50,very thick,fill=blue!20,minimum size=5mm]
\node[state] (x_0) {$p_0$};
\draw (0,1) node {$T_S$};
\node[state] (x_1) [below left of=x_0] {$p_1$};
\node[state] (x_2) [below right of=x_0] {$p_2$};
\node[state] (x_3) [below right of=x_1] {$p_3$};
\node[state] (y_1) [right of=x_2] {$q_1$};
\node[state] (y_0) [above right of=y_1] {$q_0$};
\draw (4.8,1) node {$T_P$};
\node[state] (y_2) [below right of=y_0] {$q_2$};
\node[state] (y_3) [below right of=y_1] {$q_3$};
\node[state] (z_1) [right of=y_2] {$r_1$};
\node[state] (z_0) [above right of=z_1] {$r_0$};
\draw (9.6,1) node {$T_S\parallel T_P$};
\node[state] (z_2) [below right of=z_0] {$r_2$};
\node[state] (z_3) [below right of=z_1] {$r_3$};
\path[->] (x_0) edge node {$a$} (x_1)
edge node [swap] {$a$} (x_2)
(x_1) edge node {$b$} (x_3)
(x_2) edge node [swap] {$c$} (x_3)
(y_0) edge node {$a$} (y_1)
edge node [swap] {$b$} (y_2)
(y_1) edge [bend right] node {$b$} (y_3)
          edge [bend left] node {$c$} (y_3)
(y_2) edge node [swap] {$a$} (y_3)
(z_0) edge node {$a$} (z_1)
edge node [swap] {$a$} (z_2)
(z_1) edge node {$b$} (z_3)
(z_2) edge node [swap] {$c$} (z_3);
\end{tikzpicture}
\caption{From left to right we have the transition systems modeling the specification, plant and closed-loop system. The closed-loop system is represented without states that are not reachable from the initial state.}
\label{TS}
\end{figure}
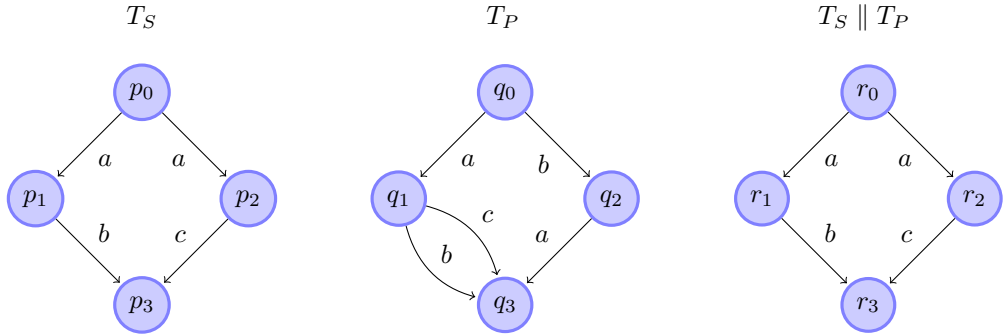

As a simple example consider the transition systems displayed in Figure~\ref{TS}. 

It is not difficult to see that the relation $R=\{(p_0,q_0),(p_1,q_1),(p_2,q_1),(p_3,q_3)\}$ is a simulation relation from $T_S$ to $T_P$. According to Corollary~\ref{Corol} there exists a controller $T_C$ making $T_C\parallel T_P$ bisimilar to $T_S$. Moreover, we know from Theorem~\ref{ExistBisimContI} that we can use $T_C=T_S$. Computing $T_S\parallel T_P$ we obtain the transition system on the right of Figure~\ref{TS} which is equal to $T_S$ and, in particular, bisimilar. A this example illustrates, even though the plant is required to be deterministic, the specification can be nondeterministic.

\subsubsection{Control systems}
In the context of control systems there are several conditions ensuring $\Po$-faithfulness of a morphism $\Sigma_1\rTo^{\sigma_{1a}}\Sigma_A$. We shall only mention the following two that will be used when discussing feedback interconnections:

\begin{enumerate}
\item the map $\sigma_{1a}=(\sigma_{1aM},\sigma_{1aU}):M_1\times U_1\to M_A\times U_A$ is injective;
\item the system $\Sigma_P$ is observable with respect to the map $\sigma_{1aU}:M\times U\to U_A$.
\end{enumerate}

Both of these assumptions guarantee that we can uniquely recover the (state) trajectory and input curve defined by a morphism $\Sigma\rTo^f \Sigma_1$ with $\Sigma$ in $\Po$ from $\Sigma\rTo^f \Sigma_1\rTo^{\sigma_{1a}} \Sigma_A$. In particular, the feedback interconnection presented in Definition~\ref{Def:ControlInterconnectI} always satisfies the first assumption.  With these considerations in place we can state the following corollary to Theorem~\ref{ExistBisimContI}.

\begin{corollary}
\label{Control}
Let $\Sigma_P=(U_P,M_P,F_P)$ and $\Sigma_S$ be control systems. The following hold:
\begin{enumerate}
\item There exists a smooth feedback control law $k:M_P\times U_C\to U_P$ making the feedback composition between $\Sigma_P$ and $k$ bisimilar to $\Sigma_S$ iff there exists a morphism $\Sigma_S\rTo^f\Sigma_P$ such that $f_M$ is a diffeomorphism.

\item Assume that $\Sigma_P$ is observable with respect to $\sigma_{1aU}$. There exists a control system $\Sigma_C$ making the feedback composition between $\Sigma_C$ and $\Sigma_P$, with interconnection maps $\phi_p$ and $\phi_s$, bisimilar to $\Sigma_S$ iff there exists a simulation relation $R$ from $\Sigma_S$ to $\Sigma_P$ satisfying:
$$(x_s,x_p)\in R\quad\implies\quad \big(F_s(x_s,(\phi_p(x_p),v_s)),F_p(x_p,(\phi_s(x_s),v_p))\big)\in TR$$
\end{enumerate}
\end{corollary}

Corollary~\ref{Control} is a straightforward instantiation of Theorem~\ref{ExistBisimContI} which nevertheless completly characterizes the solution to the controller synthesis problem using the feedback interconnections in Definitions~\ref{Def:ControlInterconnectI} and~\ref{Def:ControlInterconnectII}. These are novel results that had not been reported in the literature before. Moreover, when $F$ is control affine, $\Po$-openness of morphisms can be cheked by using the differential geometric characterizations developed in~\cite{BisimTAC,BisimSCL}.

As a simple illustration of Corollary~\ref{Control} consider the control system $\Sigma_P$ defined by:
\begin{equation}
\label{Plant}
\dot{x}=u
\end{equation}
with $x,u\in\R$, and consider also the control system $\Sigma_S$:
\begin{eqnarray}
\label{SpecI}
\dot{y}_1 & = & y_2\\
\label{SpecII}
\dot{y}_2 & = & v
\end{eqnarray}
with, $y_1,y_2,v\in\R$. Assume now that we want to construct a controller rendering control system $\Sigma_P$ bisimilar to $\Sigma_S$. We first construct the morphism $f=(f_M,f_U):\R^2\times\R\to \R\times\R$ from $\Sigma_S$ to $\Sigma_P$ by defining $f_M(y_1,y_2)=y_1$ and $f_U((y_1,y_2),v)=y_2$. The graph $R$ of $f_M$:
$$R=\left\{ ((y_1,y_2),x)\in\R^2\times\R\,\,\vert\,\, y_1=x\right\}$$
 is thus a simulation relation from the specification to the plant. This can be seen by constructing the diagram $\Sigma_S\lTo^{1_{\Sigma_S}} \Sigma_S\rTo^f \Sigma_P$ with $1_{\Sigma_S}$ the identity morphism on $\Sigma_S$ which is clearly $\Po$-open. The relation $TR$ is characterized by the equality $\dot{y}_1=\dot{x}$ or equivalently by $y_2=u$. We can thus define: 
 $$\phi_s(y)=y_2,\quad \phi_p(x)=*,\quad V_s=\{*\},\quad U_s=\R,\quad V_p=\R,\quad U_p=\{*\}$$ 
 in order to conclude that for every $(y,x)\in R$ we have $(F_s(y,v),F_p(x,\phi_s(y)))\in TR$. It now follows from observability of the plant with respect to $\sigma_{1aU}$ and from Corollary~\ref{Control} the existence of the desired controller. From Theorem~\ref{ExistBisimContI} we know that we can use the specification as the controller which in this case results in the closed loop system:
 \begin{eqnarray}
 \dot{z}_1 & = & z_2\\
 \dot{z}_1 & = & w\\
 \dot{z}_3 & = & z_2
 \end{eqnarray}
 in which we relabeled the states according to $y_1\leftrightarrow z_1$, $y_2\leftrightarrow z_2$, $y_3\leftrightarrow x$ and the input according to $v\leftrightarrow w$. The closed-loop system is easily shown to be bisimilar to the specification through the bisimulation relation:
\begin{equation}
\left\{\big((z_1,z_2,z_3),(y_1,y_2)\big)\in\R^3\times\R^2\,\,\vert\,\, z_3=y_1\land z_2=y_2\right\}
\end{equation}

\section{Further examples and discussion}
\label{Sec:Further}
\subsection{Discrete-event systems}
\label{SSec:Discrete}
In the context of supervisory control of  discrete-event systems~\cite{DEDSBook,CassandrasBook} labels are usually divided into controllable and uncontrollable. Controllable labels model transitions that can be disabled by the controller while uncontrollable labels describe the influence of the environment which is beyond the influence of control. In this setting, Theorem~\ref{ExistBisimContI} needs to be extended by adding one additional condition requiring the controller not to interfere with uncontrollable labels. Remarkably, this condition can be still be expressed in the context of open maps by suitable defining control paths and environment paths as done in~\cite{OpenTab04}. However, under the presence of uncontrollable labels, bisimulation loses some of its relevance as it fails to distinguish between controllable and uncontrollable labels. One then has to resort to alternatining bisimulation~\cite{Alternating} and, as was shown in~\cite{OpenTab04}, the framework used in this paper can still be used to prove a variant of Theorem~\ref{ExistBisimContI} that suitable takes into account uncontrollable labels. We refer the interested reader to~\cite{OpenTab04} since summarizing those results here would require us to consider the more sophisticated notion of alternating bisimulation that goes beyond the scope of this paper.

The results in Section~\ref{SSec:Existence-Examples} relied on the determinism assumption. When this assumption fails the controller synthesis problem is still solvable as shown in~\cite{NDetBisimEq}. In this more general setting the specification can no longer be used as a controller and this causes an exponential blow-up in time complexity. The exponential nature of the solution is a direct consequence of the absence of $\Po$-faithfulness since from the path in the mediating object one cannot uniquely determine the corresponding path in the plant. One is then forced to sift through all sets of possible paths in the plant corresponding to a path in the mediating object as done in~\cite{NDetBisimEq}.

A different version of Problem~\ref{CSBE}, in which bisimulation equivalence is replaced by language equivalence, has been thoroughly investigated since the pioneering work of Ramadage and Wonham~\cite{RamWonI,RamWonII}. Since for deterministic transition systems, language equivalence is equivalent to bisimulation equivalence, many of the existing results can also be obtained through a variant of Theorem~\ref{ExistBisimContI} in~\cite{OpenTab04} which distinguishes between controllable and uncontrollable labels.

\subsection{Behavioral systems}
\label{SSec:Behavioral}
In the behavioral setting~\cite{Behaviors} one considers a time set $T$, usually $\R$ or $\N$, and one describes a system $\mathcal{X}$ as a subset $\mathcal{X}\subseteq X^T$ for some set $X$. An element $\mathbf{x}\in \mathcal{X}$ is a behavior for the variable $x$ and $\mathcal{X}$ is described by the collection of all possible behaviors that $x$ may assume. Requiring all the behaviors to be defined on the same time set $T$ is restrictive since examples of nonlinear systems abound for which trajectories are only defined for sufficiently small time. We will thus take a more liberal view of a behavioral system $\mathcal{X}$ by regarding it as a subset $\mathcal{X}\subseteq \coprod_{I\in \mathbf{I}} X^I$ where $\mathbf{I}$ is the set of all intervals of the form\footnote{A similar onstruction can be performed for the discrete time case.} $]-a,b[$ with $a,b>0$.

A category of behavioral systems can be obtained by letting systems of the form $\mathcal{X}_1\subseteq \coprod_{I\in \mathbf{I}} X_1^I$ and \mbox{$\mathcal{X}_2\subseteq \coprod_{I\in \mathbf{I}} X_2^I$} be objects and by defining morphisms $\mathcal{X}_1\rTo^f\mathcal{X}_2$ as maps $f:X_1\to X_2$ taking behaviors of $\mathcal{X}_1$ into behaviors of $\mathcal{X}_2$, that is, such that for every $\mathbf{x}_1\in\mathcal{X}_1$ we have $f\circ\mathbf{x}_1\in \mathcal{X}_2$. In the behavioral setting, the composition of $\mathcal{X}_1\subseteq  \coprod_{I\in \mathbf{I}} (X_1\times Y)^I$ with $\mathcal{X}_2\subseteq  \coprod_{I\in \mathbf{I}} (X_2\times Y)^I$ through the shared variable $y\in Y$ is defined by:
$$\mathcal{X}_1\parallel_y\mathcal{X}_2=\Big\{(\mathbf{x}_1,\mathbf{x}_2)\in \coprod_{I\in\mathbf{I}} (X_1\times X_2)^I\,\,\vert\,\,\exists\, \mathbf{y}\in \coprod_{I\in\mathbf{I}} Y^I\,\,:\,\, (\mathbf{x}_1,\mathbf{y})\in \mathcal{X}_1\text{ and }(\mathbf{x}_2,\mathbf{y})\in \mathcal{X}_2\Big\}$$

This composition can also be described by a pullback. To do so we consider the system $\mathcal{A}= \coprod_{I\in\mathbf{I}} Y^I$ and the morphisms $\mathcal{X}_1\rTo^{x_{1a}}\mathcal{A}$ and $\mathcal{X}_2\rTo^{x_{2a}}\mathcal{A}$ defined by $x_{1a}(x_1,y)=y$ and $x_{2a}(x_2,y)=y$. It is not difficult to see that the pullback of $\mathcal{X}\rTo^{x_{1a}}\mathcal{A}$ and $\mathcal{X}\rTo^{x_{2a}}\mathcal{A}$ is precisely $\mathcal{X}_1\parallel_y\mathcal{X}_2$.

Since binary products exist we conclude that Theorem~\ref{ExistBisimContI} is applicable in the Behavioral context. Existence of a controller is then characterized by the existence  of a simulation relation from the specification to the plant inducing a commutative diagram such as~(\ref{SpanI}).


\subsection{Hybrid systems}
The controller synthesis problem for hybrid systems can also be solved under the proposed framework. We shall only present a brief discussion since it would take too much space to formalize all the necessary concepts. The interested reader can find such formalization in~\cite{BisimTCS} where it is shown how hybrid systems can be made into a category and how bisimulation for hybrid systems can also be described through the open maps formalism. Moreover, it is also shown in \cite{BisimTCS} that the category of hybrid systems has binary pullbacks. As expected, Theorem~\ref{ExistBisimContI} instantiated in this category implies that the controller synthesis problem is solvable when there exists a simulation relation from the specification to the plant.

\subsection{Other classes}
In the literature one can find several models of systems that have not been explicitly considered in this paper such as abstract state systems~\cite{PSDB05}, general systems in the behavioral setting~\cite{vdS03} and general flow systems~\cite{DT07} among many others. Provided that the corresponding categories have binary products the results presented in this paper also bring considerable insight into these specific classes of systems.

\subsection{Discussion}
The controller synthesis problem for bisimulation equivalence admits a very intuitive and simple solution that is valid across a wide range of systems: \emph{a controller exists iff the plant simulates the specification}. The simplicity of this statement is a consequence of the categorical approach taken in this paper that distilled the essence of the problem and lead to a solution bringing considerable insight into the concrete classes of systems to which it was applied. The proposed solution also points to the need of developping computational efficient methods to determine the existence of simulation relations between several classes of systems including control and hybrid systems.
%

\bibliographystyle{alpha}
\bibliography{BisimEq}

\end{document}